\documentclass[a4paper,9pt]{article}
\usepackage{amsmath, amsfonts, amssymb}
\usepackage[english, russian]{babel}
\usepackage[cp1251]{inputenc}

\usepackage{graphicx}
\begin{document}
\sloppy
\begin{center}
\textbf{New decomposition formulas associated with the Lauricella multivariable hypergeometric functions }\\[5pt]
% Авторы
\textbf{Ergashev T.G.\\}
% обязательное поле!
\medskip
{ Institute of Mathematics, Uzbek
Academy of Sciences,  Tashkent, Uzbekistan. \\

{\verb ergashev.tukhtasin@gmail.com }}\\
\end{center}

\begin{quote} % Аннотация на английском

Decomposition formulas associated with the Lauricella multivariable hypergeometric functions were known, however, due to the recurrence of those formulas, additional difficulties may
arise in the applications. Further study of the properties
of the famous expansion formulas showed that it can be reduced to a
more convenient form. In addition, this paper contains applications of new expansion formulas to the solving of boundary value problems for a multidimensional elliptic equation with several singular coefficients.

\end{quote}

\section{Introduction}

A great interest in the theory of multiple hypergeometric functions is motivated  essentially by the fact that the solutions of many applied problems involving (for example) partial differential equations are obtainable with the help of such hypergeometric functions (see, for details, \cite[p.47 et seq. Section 1.7]{SK}; see also the works \cite{{Opps},{Padma}} and the references cited therein). For instance, the energy absorbed by some nonferromagnetic conductor sphere included in an  internal mgnetic field can be calculated with the help of such functions \cite{Loh}. Hypergeometric functions of  several variables are used in physical and quantum chemical applications as well \cite{{Niuk},{Opps}} . Especially, many problems in gas dynamics lead to solutions of degenerate second-order partial differential equations which are then solvable in terms of multiple hypergeometric functions. Among examples, we can cite the problem of adiabatic flat-parallel gas flow without whirlwind, the flow problem of supersonic current from vessel with flat walls, and a number of other problems connected with gas flow \cite{Frankl}.

We note that Riemann's functions, Green's functions and the fundamental solutions of the degenerate second-order partial differential  equations are expressible by means of hypergeometric functions of several variables \cite{{BN1},{BN2},{BN3},{Erg},{EH},{H},{HK},{M},{U}, {UrK},{Wt}}. In investigation of the boundary-value problems for these partial differential equations, we need decompositions for hypergeometric functions of several variables in terms of simpler hypergeometric functions of (for example) the Gauss and Appell types.

The familiar operator method of Burcnall and Chaundy \cite{{BC1},{BC2},{Chaundy}} has been used by them rather extensively for finding  decomposition formulas for hypergeometric functions of two variables in terms of the classical Gauss hypergeometric function of one variable.

Following the works \cite{{BC1},{BC2}}, Hasanov and Srivastava \cite{{HS6},{HS7}} introduced operators generalizing the Burcnall-Chaundy operators and found expansion formulas for many triple hypergeometric functions which were successfully applied to the solving the boundary-value problems for the second order elliptic equation with three singular coefficients \cite{{Kar},{A20},{A22}}, and they proved recurrent formulas when the dimension of hypergeometric function exceeds three. However, due to the recurrence, additional difficulties may
arise in the applications of those decomposition formulas.

In this paper for the two Lauricella hypergeometric functions in several variables we prove new decomposition formulas which are free from the recurrence and applied to the solving the boundary-value problems for the multidimensional elliptic equation with several singular coefficients.

The plan of this paper is as follows. In Section 2 we briefly give some preliminary information, which will be used later.
In Section 3, we present the well-known decomposition formulas associated with the two and more dimensional Lauricella hypergeometric functions. In Section 4, we will prove new decomposition and summation formulas and in the last section 5 we will apply the obtained formulas to the solution of boundary-value problems.

\section{Preliminaries}

Below we give some formulas for Euler gamma-function, Gauss hypergeometric function, Lauricella hypergeometric functions of three and more variables, which will be used in the next sections.

Let be ${N}$ set of the natural numbers : $N=\{1,2,3,...\}.$

It is known that the Euler gamma-function $\Gamma(a)$ has property \cite[p.17, (2)]{Erd}
$$\Gamma(a+m)=\Gamma(a) (a)_m.$$
Here $(a)_m$ is a Pochhammer symbol, for which the equality
$(a)_{m+n}=(a)_m(a+m)_n $ and its particular case $(a)_{2m}=(a)_m(a+m)_m$ are true \cite[p.67,(5)]{Erd}.

A function
$$F {\left(a,b;c;x\right)} \equiv F {\left[
{{\begin{array}{*{20}c}
 {a,b} ; \hfill \\
 {c};  \hfill \\
\end{array}} x}  \right]}={\sum\limits_{i =
0}^{\infty} \frac{(a)_i(b)_i}{(c)_ii!}x^i }, \,c\neq 0,-1,-2,...$$
is known as the  Gaussian hypergeometric function and an equality
\begin{equation}
\label{sum}
F(a,b;c;1)=\frac{\Gamma(c)\Gamma(c-a-b)}{\Gamma(c-a)\Gamma(c-b)}, c\neq 0,-1,-2,..., Re(c-a-b)>0
\end{equation}
holds \cite[p.73, (73)]{Erd}. Moreover, the following autotransformation formula \cite[p.76, (22)]{Erd}
\begin{equation}
\label{auto}
F\left(a,b;c;x\right)=\left(1-x\right)^{-b}F\left(c-a,b;c;\frac{x}{x-1}\right)
\end{equation}
is valid.

Multiple Lauricella  hypergeometric functions $F_A^{(n)}$ and $F_B^{(n)}$ in $n\in N$ (real or complex) variables are defined as following (\cite{A30} and \cite[p.33]{AP})
\begin{equation*}
F_{A}^{(n)} \left( {a,b_{1} ,...,b_{n} ;c_{1} ,...,c_{n} ;x_{1} ,...,x_{n}}
\right) \equiv F_{A}^{(n)} {\left[ {{\begin{array}{*{20}c}
 {a,b_{1} ,...,b_{n} ;} \hfill \\
 {c_{1} ,...,c_{n} ;} \hfill \\
\end{array}} x_{1} ,...,x_{n}}  \right]}
\end{equation*}
\begin{equation*}
\label{eq7}
 = {\sum\limits_{m_{1} ,...m_{n} = 0}^{\infty}  {{\frac{{\left( {a}
\right)_{m_{1} + ... + m_{n}}  \left( {b_{1}}  \right)_{m_{1}}  ...\left(
{b_{n}}  \right)_{m_{n}}} } {{\left( {c_{1}}  \right)_{m_{1}}  ...\left(
{c_{n}}  \right)_{m_{n}}} } }{\frac{{x_{1}^{m_{1}}} } {{m_{1}
!}}}...{\frac{{x_{n}^{m_{n}}} } {{m_{n} !}}}}}
\end{equation*}
\[\,{\left[ {c_{k} \ne 0,-1,-2,...;\,k = \overline {1,n} ;\,{\left| {x_{1}}
\right|} + ... + {\left| {x_{n}}  \right|} < 1} \right]};\]

\begin{equation*}
F_{B}^{(n)} \left( {a_1,...,a_n,\,b_{1} ,...,b_{n} ;c;x_{1} ,...,x_{n}}
\right) \equiv F_{B}^{(n)} {\left[ {{\begin{array}{*{20}c}
 {a_1,...,a_n,\,b_{1} ,...,b_{n} ;} \hfill \\
 {c;} \hfill \\
\end{array}} x_{1} ,...,x_{n}}  \right]}
\end{equation*}
\begin{equation*}
\label{eq71}
 = \sum\limits_{m_{1} ,...m_{n} = 0}^{\infty}  {{\frac{{\left( {a_1}
\right)_{m_{1}}...\left( {a_1}
\right)_{m_{n}}  \left( {b_{1}}  \right)_{m_{1}}  ...\left(
{b_{n}}  \right)_{m_{n}}} } {{\left({c}\right)_{m_{1}+...+m_n} } }{\frac{{x_{1}^{m_{1}}} } {{m_{1}
!}}}...{\frac{{x_{n}^{m_{n}}} } {{m_{n} !}}}}}
\end{equation*}
\[\,\left[ {c \ne 0,-1,-2,... ;\,{max\{\left| {x_{1}}
\right|,...,{\left| {x_{n}}  \right|}\} < 1} }\right].\]

\section{Decomposition formulas associated with the Lauricella functions $F_A^{(n)}$ and $F_B^{(n)}$ }

For a given multiple hypergeometric function, it is useful to fund a
decomposition formula which would express the multivariable hypergeometric
function in terms of products of several simpler hypergeometric functions
involving fewer variables.

Burchnall and Chaundy \cite{{BC1},{BC2}} systematically
presented a number of expansion and decomposition formulas for some double
hypergeometric functions in series of simpler hypergeometric functions. For
example, the Appell function
\begin{equation*}
F_{2} \left( {a,b_{1} ,b_{2} ;c_{1} ,c_{2} ;x,y} \right) = {\sum\limits_{i,j
= 0}^{\infty}  {{\frac{{\left( {a} \right)_{i+j} \left( {b_{1}}
\right)_{i} \left( {b_{2}}  \right)_{j}}} {{\left( {c_{1}}  \right)_{i}
\left( {c_{2}}  \right)_{j}}} }{\frac{{x^{i}}}{{i!}}}{\frac{{y^{j}}}{{j!}}}}
}
\end{equation*}
\begin{equation*}
{\left[ {c_{1} ,c_{2} \ne 0,-1,-2,...;\,\,{\left| {x} \right|} + {\left| {y} \right|}
< 1} \right]}
\end{equation*}
\noindent
has the expansion \cite{BC1}
\[
F_{2} \left( {a,b_{1} ,b_{2} ;c_{1} ,c_{2} ;x,y} \right)\]
\begin{equation}
\label{eq25}
= {\sum\limits_{i =
0}^{\infty}  {{\frac{{\left( {a} \right)_{i} \left( {b_{1}}  \right)_{i}
\left( {b_{2}}  \right)_{i}}} {{i!\left( {c_{1}}  \right)_{i} \left( {c_{2}
} \right)_{i}}} }x^{i}y^{i}F\left( {a + i,b_{1} + i;c_{1} + i;x}
\right)F\left( {a + i,b_{2} + i;c_{2} + i;y} \right)}}.
\end{equation}

The Burchnall-Chaundy method, which is limited to functions of two variables,
is based on the following mutually inverse symbolic operators \cite{BC1}
\begin{equation}
\label{eq8}
\nabla \left( {h} \right) = {\frac{{\Gamma \left( {h} \right)\Gamma \left(
{{\rm \delta} _{1} + {\rm \delta} _{2} + h} \right)}}{{\Gamma \left( {{\rm
\delta} _{1} + h} \right)\Gamma \left( {{\rm \delta} _{2} + h} \right)}}},
\quad
\Delta \left( {h} \right) = {\frac{{\Gamma \left( {{\rm \delta} _{1} + h}
\right)\Gamma \left( {{\rm \delta} _{2} + h} \right)}}{{\Gamma \left( {h}
\right)\Gamma \left( {{\rm \delta} _{1} + {\rm \delta} _{2} + h}
\right)}}},
\end{equation}
\noindent
where ${\rm \delta} _{1} = x{\displaystyle\frac{{\partial}} {{\partial x}}}$ and ${\rm
\delta} _{2} = y{\displaystyle\frac{{\partial}} {{\partial y}}}$.

In order to generalize the operators $\nabla \left( {h} \right)$ and $\Delta
\left( {h} \right)$, defined in (\ref{eq8}), Hasanov and Srivastava \cite{{HS6}, {HS7}}
introduced the operators
\begin{equation*}
\label{eq9}
\tilde {\nabla} _{x_{1} ;x_{2} ,...,x_{n}}  \left( {h} \right) =
{\frac{{{\rm \Gamma} \left( {h} \right){\rm \Gamma} \left( {{\rm \delta
}_{1} + ... + {\rm \delta} _{n} + h} \right)}}{{{\rm \Gamma} \left( {{\rm
\delta} _{1} + h} \right){\rm \Gamma} \left( {{\rm \delta} _{2} + ... + {\rm
\delta} _{n} + h} \right)}}},
\end{equation*}
\begin{equation*}
\label{eq10}
\tilde {\Delta} _{x_{1} ;x_{2} ,...,x_{n}}  \left( {h} \right) =
{\frac{{{\rm \Gamma} \left( {{\rm \delta} _{1} + h} \right){\rm \Gamma
}\left( {{\rm \delta} _{2} + ... + {\rm \delta} _{n} + h} \right)}}{{{\rm
\Gamma} \left( {h} \right){\rm \Gamma} \left( {{\rm \delta} _{1} + ... +
{\rm \delta} _{n} + h} \right)}}},
\end{equation*}
\noindent
where ${\rm \delta} _{k} = x_{k} {\displaystyle\frac{{\partial}} {{\partial x_{k}}} }\,$ $(k=\overline{1,n}),$
with the help of which they managed to find decomposition formulas for a
whole class of hypergeometric functions in several variables.

Following the works  \cite{{BC1},{BC2}} Hasanov and Srivastava
\cite{HS6} found following decomposition formulas for the
Lauricella functions of three variables
\begin{equation}\label{e26}\begin{array}{l}F^{(3)}_A \left( {{
 {a,b_1,b_2,b_3;}{c_1,c_2,c_3;}} x_1,x_2,x_3}  \right)
={\sum\limits_{i,j,k=0}^\infty\displaystyle\frac{(a)_{i+j+k}(b_1)_{j+k}(b_2)_{i+k}(b_3)_{i+j}}{i!j!k!(c_1)_{j+k}(c_2)_{i+k}(c_3)_{i+j}}}x_1^{j+k}x_2^{i+k}x_3^{i+j}
\\
\\
\cdot
F\left( {{{a+j+k,b_1+j+k;} {c_1+j+k;} } x_1}  \right)F\left( {{{a+i+j+k,b_2+i+k;}{c_2+i+k;}} x_2}  \right)\\
\\
\cdot F\left( {{{a+i+j+k,b_3+i+j;}{c_3+i+j;}} x_3}  \right), \end{array}
\end{equation}
\begin{equation*}
\label{eq66}
\begin{array}{l}
 F_{B}^{(3)} \left( {a_{1} ,a_{2} ,a_{3} ;b_{1} ,b_{2} ,b_{3} ;c;x_{1}
,x_{2} ,x_{3}}  \right) \\
\\
 = {\sum\limits_{i,j,k = 0}^{\infty}  {{\displaystyle\frac{{\left( { - 1} \right)^{i + j
+ k}\left( {a_{1}}  \right)_{j + k} \left( {b_{1}}  \right)_{j + k} \left(
{a_{2}}  \right)_{i + k} \left( {b_{2}}  \right)_{i + k} \left( {a_{3}}
\right)_{i + j} \left( {b_{3}}  \right)_{i + j}}} {{\left( {c - 1 + j + k}
\right)_{j + k} \left( {c - 1 + 2\left( {j + k} \right) + i} \right)_{i}
\left( {c} \right)_{2\left( {i + j + k} \right)} i!j!k!}}}}}x_{1}^{j + k} x_{2}^{i + k} x_{3}^{i + j}  \\
\\
\cdot  F\left( {a_{1} + j + k,b_{1} + j + k;c + 2\left( {j + k} \right);x_{1}}
\right)F\left( {a_{2} + i + k,b_{2} + i + k;c + 2\left( {i + j + k} \right);x_{2}
} \right) \\
\\
\cdot F\left( {a_{3} + i + j,b_{3} + i + j;c + 2\left( {i + j + k}
\right);x_{3}}  \right) \\
 \end{array}
\end{equation*}
and they proved that for all $n\in N\backslash\{1\}$ are true the
recurrence formulas \cite {HS7}
\begin{equation}
\label{e27} \begin{array}{l} F_{A}^{(n)}\left(
{{{a,b_1,...,b_n;}{c_1,...,c_n;}} x_1,...,x_n}  \right) \\
\\
 = {\sum\limits_{m_{2} ,...,m_{n} = 0}^{\infty}  {{\displaystyle\frac{{(a)_{m_{2} + \cdot
\cdot \cdot + m_{n}}  (b_{1} )_{m_{2} + \cdot \cdot \cdot + m_{n}}
(b_{2} )_{m_{2}}  \cdot \cdot \cdot (b_{n} )_{m_{n}}} } {{m_{2} !
\cdot \cdot \cdot m_{n} !(c_{1} )_{m_{2} + \cdot \cdot \cdot +
m_{n}}  (c_{2} )_{m_{2}}  \cdot \cdot \cdot (c_{n} )_{m_{n}}} }
}}} x_{1}^{m_{2} + \cdot \cdot \cdot + m_{n}} x_{2}^{m_{2}}
\cdot \cdot \cdot x_{n}^{m_{n}}  \\
\\
\cdot x_{1}^{m_{2} + \cdot \cdot \cdot + m_{n}} F\left( {{{a + m_{2} + \cdot \cdot \cdot + m_{n},b_{1} + m_{2} +
\cdot \cdot \cdot
+ m_{n};}
 {c_{1} + m_{2} +
\cdot \cdot \cdot
+ m_{n};}
} x_1}  \right) \\
\\
 \cdot F_{A}^{(n - 1)} \left( {{ {a + m_{2} + \cdot \cdot \cdot + m_{n} ,b_{2}
+ m_{2} ,...,b_{n} + m_{n} ;}
c_{2} + m_{2} ,....,c_{n} + m_{n}
; } x_{2} ,...,x_{n}}  \right),
\end{array}
\end{equation}

\begin{equation}
\label{e277}
\begin{array}{l}
 F_{B}^{(n)} \left( {a_{1} ,...,a_{n} ;b_{1} ,...,b_{n} ;c;x_{1} ,...,x_{n}
} \right) \\
\\
 = {\sum\limits_{k_{2} ,...,k_{n} = 0}^{\infty}  {{\displaystyle\frac{{\left( { - 1}
\right)^{k_{2} + ... + k_{n}} \left( {a_{1}}  \right)_{k_{2} + ... + k_{n}}
\left( {b_{1}}  \right)_{k_{2} + ... + k_{n}}  {\prod\nolimits_{j = 2}^{n}
{{\left[ {\left( {a_{j}}  \right)_{k_{j}}  \left( {b_{j}}  \right)_{k_{j}}
} \right]}}}} }{{\left( {c - 1 + k_{2} + ... + k_{n}}  \right)_{k_{2} + ...
+ k_{n}}  \left( {c} \right)_{2\left( {k_{2} + ... + k_{n}}  \right)} k_{2}
!...k_{n} !}}}}} \\
\\
\cdot  x_{1}^{k_{2} + ... + k_{n}}  x_{2}^{k_{2}}  ...x_{n}^{k_{n}
} F\left( {a_{1} + k_{2} + ... + k_{n} ,b_{1} + k_{2} + ... + k_{n} ;c +
2\left( {k_{2} + ... + k_{n}}  \right);x_{1}}  \right) \\
\\
\cdot F_{B}^{(n - 1)} \left( {a_{2} + k_{2} ,...,a_{n} + k_{n} ,b_{2} + k_{2}
,...,b_{n} + k_{n} ;c + 2\left( {k_{2} + ... + k_{n}}  \right);x_{2}
,...,x_{n}}  \right) .\\
 \end{array}
\end{equation}

However, due to the recurrence of formula (\ref{e27}) and (\ref{e277}), additional difficulties may
arise in the applications of this expansion. Further study of the properties
of the Lauricella functions $F_A^{(n)}$ and $F_B^{(n)}$  showed that formulas (\ref{e27}) and (\ref{e277}) can be reduced to a
more convenient forms.

\bigskip

\section{New decomposition formulas associated with the Lauricella functions $F_A^{(n)}$ and $F_B^{(n)}$ }

Before proceeding to the presentation of the main result of this article, we introduce the notations
\begin{equation}
\label{e4111}
A(k,n)=\sum\limits_{i=2}^{k+1}\sum\limits_{j=i}^n m_{i,j},\,\,B(k,n)=\sum\limits_{i=2}^{k}m_{i,k}+\sum\limits_{i=k+1}^n m_{k+1,i},
\end{equation}
where  $m_{i,j}\in \rm{N}\cap\{0\} \left(2\leq i\leq j \leq n\right).$

It should be noted here that the sum ${ {B(2,n)+B(3,n)+...+B(n,n)}} $
has the parity property, which plays an important role in the calculation of
the some values of hypergeometric functions. In fact, by virtue of equality
\begin{equation*}
{\sum\limits_{k = 2}^{n} {{\sum\limits_{i = 2}^{k} {m_{i,k}}} } }  =
{\sum\limits_{k = 1}^{n - 1} {{\sum\limits_{i = k + 1}^{n} {m_{k + 1,i}}} }
}
\end{equation*}
\noindent
we obtain
\begin{equation}
\label{eq1444}
{\sum\limits_{k = 1}^{n} {B(k,n)}}  = 2{\sum\limits_{k = 2}^{n}
{{\sum\limits_{i = 2}^{k} {m_{i,k}}} } }  = 2{\sum\limits_{k = 1}^{n - 1}
{{\sum\limits_{i = k + 1}^{n} {m_{k + 1,i}}} } } .
\end{equation}

We present other simple properties of the functions $A\left(
{k,n} \right)$ and $B\left( {k,n} \right)$:
\begin{equation}
\label{eq15555}
A\left( {n + 1,n + 1} \right) - B\left( {n + 1,n + 1} \right) = A\left(
{n,n} \right),
\end{equation}
\begin{equation}
\label{eq16666}
A\left( {k + 1,k + 1} \right) - B\left( {k + 1,k + 1} \right) = A\left(
{k,n} \right) - B\left( {k,n} \right) + m_{2,n + 1} + ... + m_{k,n + 1} .
\end{equation}

Those properties are easily proved if we proceed from the
definitions of functions $A\left( {k,n} \right)$ and $B\left( {k,n}
\right)$.

\textbf{Lemma 1}. The following decomposition formulas hold true at $n
\in {\rm N}$
\begin{equation*}
F_{A}^{(n)} \left( {a,b_{1} ,b_{2} ,....,b_{n} ;c_{1} ,c_{2} ,....,c_{n}
;x_{1} ,...,x_{n}}  \right)
\end{equation*}
\begin{equation}
\label{eq1222} 
 = {\sum\limits_{{\mathop {m_{i,j} = 0}\limits_{(2 \le i \le j \le n)}
}}^{\infty}  {{\frac{{(a)_{A(n,n)}}} {{m_{ij} !}}}}} \prod\limits_{k =
1}^{n} [ \frac{{(b_{k} )_{B(k,n)}}}{(c_{k})_{B(k,n)}}
 x_{k}^{B(k,n)}
F( a + A(k,n), b_{k} + B(k,n);c_{k} + B(k,n);x_{k})],
\end{equation}

\begin{equation*}
 F_{B}^{(n)} \left( {a_{1} ,...,a_{n} ,b_{1} ,...,b_{n} ;c;x_{1} ,...,x_{n}
} \right)
\end{equation*}

\begin{equation*}
 = {\sum\limits_{{\mathop {m_{i,j} = 0}\limits_{(2 \le i \le j \le n)}
}}^{\infty}  {{\displaystyle\frac{{\left( { - 1} \right)^{A\left( {n,n}
\right)}}}{{\left( {c} \right)_{2A\left( {n,n} \right)} m_{ij} !}}}}
}\prod\limits_{k = 1}^{n}  [{\displaystyle\frac{{\left( {a_{k}}
\right)_{B\left( {k,n} \right)} \left( {b_{k}}  \right)_{B\left( {k,n}
\right)} \left( {c - 1} \right)_{A(k,n) - A(k - 1,n)}}} {{\left( {c - 1}
\right)_{2A(k,n) - 2A(k - 1,n)}}} }
\end{equation*}
\begin{equation}
\label{eq1111} \cdot x_{k}^{B\left( {k,n} \right)} F\left(
{a_{k} + B\left( {k,n} \right),b_{k} + B\left( {k,n} \right);c + 2A\left(
{k,n} \right);x_{k}}  \right)].
\end{equation}

\textbf{Proof}. We carry out the proof by the method mathematical
induction. First, we prove the validity of the equality (\ref{eq1222}).

For clarity of the course of the proof, we introduce the notations
$$
N_{l} (k,n)
= {\sum\limits_{i = l}^{k + 1} {{\sum\limits_{j = i}^{n}
{m_{i,j}}} } } , \, M_{l} (k,n) = {\sum\limits_{i = l}^{k} {m_{i,k} +}}
{\sum\limits_{i = k + 1}^{n} {m_{k + 1,i}}}, l \in N.
$$

It's obvious that
$$
N_{2} (k,n) = A(k,n),\, M_{2} (k,n)=B(k,n).
$$

So we have to prove the fairness of equality
\begin{equation}
\label{e28}\begin{array}{l} F_{A}^{(n)} {\left[
{{\begin{array}{*{20}c}{a,b_{1}
,....,b_{n} ;} \hfill \\ {c_{1} ,....,c_{n} ;} \hfill \\
\end{array}} x_{1} ,...,x_{n}}  \right]}
 = {\sum\limits_{{\mathop {m_{i,j} = 0}\limits_{(2 \le i \le j \le n)}
}}^{\infty}  {{\displaystyle\frac{{(a)_{N_{2} (n,n)}}} {{ m_{i,j}!}}}}} \hfill\\
\\
 \cdot{\prod\limits_{k = 1}^{n} {{ {{\displaystyle\frac{{(b_{k} )_{M_{2} (k,n)}
}}{{(c_{k} )_{M_{2} (k,n)}}} x_{k}^{M_{2} (k,n)}
F\left[{\begin{array}{*{20}c} {a + N_{2} (k,n),b_{k} + M_{2}
(k,n);} \hfill \\ c_{k} + M_{2} (k,n);\hfill \\ \end{array} x_{k}}
\right]} } }}}.
\end{array}
\end{equation}

In the case $n=1$ the equality (\ref{e28}) is obvious.

Let $n = 2$. Since $M_2(1,2)=M_2(2,2)=N_2(1,2)=N_2(2,2)=m_{2,2}:= i,$
we obtain the formula (\ref{eq25}).

For the sake of interest, we will check the formula (\ref{e28}) in
yet another value of $n$.

Let $n=3.$ In this case
$$M_2(1,3)=m_{2,2}+m_{2,3},\,\, M_2(2,3)=m_{2,2}+m_{3,3},\,\, M_2(3,3)=m_{2,3}+m_{3,3},$$
$$N_2(1,3)=m_{2,2}+m_{2,3},\,\, N_2(2,3)= N_2(3,3)=m_{2,2}+m_{2,3}+m_{3,3}.$$
For brevity, making the substitutions
$m_{2,2}:=i,\,\,m_{2,3}:=j,\,\,m_{3,3}:=k$, we obtain the formula
(\ref{e26}).

So the formula (\ref{e28}), that is formula (\ref{eq1222}), works for $n=1,$ $n=2$ and $n=3$.

Now we assume that for $n = s$ equality (\ref{e28}) holds; that
is, that
\begin{equation}
\label{e29}
\begin{array}{l}
 F_{A}^{(s)} {\left[ \begin{array}{*{20}c}{a,b_{1} ,....,b_{s} ;}\hfill\\ c_{1} ,....,c_{s} ; \hfill\\ \end{array} x_{1}
,...,x_{s}  \right]}
 = {\sum\limits_{{\mathop {m_{i,j} = 0}\limits_{(2 \le i \le j \le s)}
}}^{\infty}  {{\displaystyle\frac{{(a)_{N_{2} (s,s)}}} {m_{ij}! }}}}  \\
\\
 \cdot {\prod\limits_{k = 1}^{s} {{{{\displaystyle\frac{{(b_{k} )_{M_{2} (k,s)}
}}{{(c_{k} )_{M_{2} (k,s)}}} }x_{k}^{M_{2} (k,s)} F\left[
\begin{array}{*{20}c}{a + N_{2} (k,s),b_{k} + M_{2} (k,s);}\hfill\\ c_{k} +
M_{2} (k,s);\hfill\\ \end{array}  x_{k} \right]} }}} .
\\
 \end{array}
\end{equation}

Let $n=s+1.$ We prove that  following formula
\begin{equation}
\label{e210}
\begin{array}{l}
 F_{A}^{(s + 1)} {\left[ \begin{array}{*{20}c}{a,b_{1} ,....,b_{s+1} ;}\hfill\\ c_{1} ,....,c_{s+1} ; \hfill\\ \end{array} x_{1}
,...,x_{s+1}  \right]}
 = {\sum\limits_{{\mathop {m_{i,j} = 0}\limits_{(2 \le i \le j \le s + 1)}
}}^{\infty}  {{\displaystyle\frac{{(a)_{N_{2} (s + 1,s + 1)}}} {m_{ij}! }}}}  \\
\\
 \cdot {\prod\limits_{k = 1}^{s + 1} {{ {{\displaystyle\frac{{(b_{k} )_{M_{2} (k,s + 1)}
}}{{(c_{k} )_{M_{2} (k,s + 1)}}} }x_{k}^{M_{2} (k,s + 1)} F\left[
{{\begin{array}{*{20}c}
 {a + N_{2} (k,s + 1),b_{k} + M_{2} (k,s + 1);} \hfill \\
 {c_{k} + M_{2} (k,s + 1);} \hfill \\
\end{array}} x_{k}}  \right]} }}}  \\
 \end{array}
\end{equation}
is valid.

We write the Hasanov-Srivastava's formula (\ref{e27}) in the form

\begin{equation}
\label{e211}
\begin{array}{l}
 F_{A}^{(s + 1)} {\left[ {{\begin{array}{*{20}c}{a,b_{1} ,....,b_{s + 1} ;}\hfill \\ { c_{1} ,....,c_{s +
 1};}\hfill \\
 \end{array}}
x_{1} ,...,x_{s + 1}}  \right]} \\
\\
 = {\sum\limits_{m_{2,2} ,...,m_{2,s + 1} = 0}^{\infty}  {{\displaystyle\frac{{(a)_{N_2(1,s+1)}  (b_{1} )_{M_2(1,s+1)}  (b_{2} )_{m_{2,2}}  \cdot \cdot \cdot (b_{s + 1}
)_{m_{2,s + 1}}} } {{m_{2,2} ! \cdot \cdot \cdot m_{2,s + 1}
!(c_{1} )_{M_2(1,s+1)}  (c_{2} )_{m_{2,2}} \cdot \cdot \cdot (c_{s
+ 1})_{m_{2,s + 1}}} } }}}   \\
\\
 \cdot x_{1}^{M_2(1,s+1)}  x_{2}^{m_{2,2}} \cdot
\cdot \cdot x_{s + 1}^{m_{2,s + 1}} F\left[\begin{array}{*{20}c}{a + N_2(1,s+1) ,b_{1} + M_2(1,s+1) ;} \hfill\\ c_{1} +M_2(1,s+1);\hfill\\ \end{array} x_{1} \right] \\
\\
 \cdot F_{A}^{(s)} {\left[ {{\begin{array}{*{20}c}
 {a + N_2(1,s+1) ,b_{2} + m_{2,2} ,...,b_{s + 1} +
m_{2,s + 1} ;} \hfill \\
 {c_{2} + m_{2,2} ,....,c_{s + 1} + m_{2,s + 1}};  \hfill \\
\end{array}} x_{2} ,...,x_{s + 1}}  \right]}. \\
 \end{array}
\end{equation}

By virtue of  the formula (\ref{e29}) we have
\begin{equation}
\label{e212}
\begin{array}{l}
 F_{A}^{(s)} \left[ \begin{array}{*{20}c}{a + N_2(1,s+1) ,b_{2} + m_{2,2} ,...,b_{s +
1} + m_{2,s + 1} ;}\hfill\\ c_{2} + m_{2,2} ,...,c_{s + 1} +
m_{2,s + 1} ;\hfill\\ \end{array} x_{2} ,...,x_{s +
1}  \right] \\
\\
 = {\sum\limits_{{\mathop {m_{i,j} = 0}\limits_{(3 \le i \le j \le s + 1)}
}}^{\infty}  {{\displaystyle\frac{{\left( {a +N_2(1,s+1)} \right)_{N_{3} (s +
1,s + 1)}}} {m_{ij}!}}}} \prod\limits_{k = 2}^{s + 1}
{\displaystyle\frac{{(b_{k} + m_{2,k} )_{M_{3} (k,s + 1)}}} {{(c_{k} + m_{2,k}
)_{M_{3} (k,s + 1)}}}
}x_{k}^{M_{3} (k,s + 1)}  \\
\\
 \cdot F\left[ {{\begin{array}{*{20}c}
 {a + N_2(1,s+1) + N_{3} (k,s + 1),b_{k} + m_{2,k} + M_{3} (k,s +
1);} \hfill \\
 {c_{k} + m_{2,k} + M_{3} (k,s + 1);} \hfill \\
\end{array}} x_{k}}  \right]  . \\
 \end{array}
\end{equation}

Substituting from (\ref{e212}) into (\ref{e211}) we obtain
\begin{equation*}
\begin{array}{l}
 F_{A}^{(s + 1)} {\left[ {a,b_{1} ,....,b_{s + 1} ;c_{1} ,....,c_{s + 1}
;x_{1} ,...,x_{s + 1}}  \right]} \\
\\
 = {\sum\limits_{{\mathop {m_{i,j} = 0}\limits_{(2 \le i \le j \le s + 1)}
}}^{\infty}  {{\displaystyle\frac{{\left( {a} \right)_{N_2(1,s+1) + N_{3} (s +
1,s + 1)}}} {m_{ij}! }}}} \prod\limits_{k = 1}^{s + 1}
{\displaystyle\frac{{(b_{k} )_{m_{2,k} + M_{3} (k,s + 1)}}} {{(c_{k} )_{m_{2,k}
+ M_{3} (k,s + 1)}}}
}x_{k}^{m_{2,k} + M_{3} (k,s + 1)}   \\
\\
\cdot F\left[ {{\begin{array}{*{20}c} {a + N_2(1,s+1) + N_{3} (k,s
+ 1),}
 {b_{k} + m_{2,k} + M_{3} (k,s + 1);}
  \hfill \\
 {c_{k} + m_{2,k} + M_{3} (k,s + 1)}; \hfill \\
\end{array}} x_{k}}  \right]. \\
 \end{array}
\end{equation*}

Further, by virtue of the following obvious equalities
$$
N_2(1,s+1) + N_{3} (k,s + 1) = N_{2} (k,s + 1),\,\,1\leq k\leq
s+1, s\in N, $$ $$ m_{2,k} + M_{3} (k,s + 1) = M_{2} (k,s + 1),
1\leq k\leq s+1, s\in N,
$$
we finally find the equality (\ref{e210}).

The  equality (\ref{eq1111}) is proved similarly as proof of the equality (\ref{eq1222}). Q.E.D.

\bigskip

\textbf{Lemma 2}. Let $a,b_{1} ,$\ldots , $b_{n} $ are real numbers with $a = 0,\,
- 1,\, - 2,...$ and $a > b_{1} + ... + b_{n}.$ Then the following summation
formulas hold true at $n \in {\rm N}$
$$
{\sum\limits_{{\mathop {m_{i,j} = 0}\limits_{(2 \le i \le j \le n)}
}}^{\infty}  {{\displaystyle\frac{{(a)_{A (n,n)}}} {{m_{ij}!} }}}} {\prod\limits_{k = 1}^{n}\left[ {{\frac{{\left( {b_{k}}  \right)_{B(k,n)}
\left( {a - b_{k}}  \right)_{A(k,n) - B(k,n)}}} {{\left( {a}
\right)_{A(k,n)}}} }}\right]}
$$
\begin{equation}
\label{eq888}
=\frac{\Gamma \left( {a - {\sum\nolimits_{k = 1}^{n} {b_{k}
}}}  \right)}{\Gamma(a)}\prod\limits_{k = 1}^{n}\left[\frac{{\Gamma(a)}}{{{{\Gamma \left( {a - b_{k}}  \right)}}} }\right],
\end{equation}
$$
 {\sum\limits_{{\mathop {m_{i,j} = 0}\limits_{(2 \le i \le j \le n)}
}}^{\infty}  {{\displaystyle\frac{{\left( { - 1} \right)^{A\left( {n,n}
\right)}}}{{\left( {a} \right)_{2A\left( {n,n} \right)} m_{ij} !}}}}
}{\prod\limits_{k = 1}^{n} {{\left[ {{\frac{{\left( {b_{k}}
\right)_{B\left( {k,n} \right)} \left( {a} \right)_{2A\left( {k,n} \right)}
\left( {a - 1} \right)_{A(k,n) - A(k - 1,n)}}} {{\left( {c - b_{k}}
\right)_{2A\left( {k,n} \right) - B\left( {k,n} \right)} \left( {a - 1}
\right)_{2A(k,n) - 2A(k - 1,n)}}} }} \right]}}}  \\
$$
\begin{equation}
\label{eq8888}
 = {\displaystyle\frac{{\Gamma \left( {a} \right)}}{{\Gamma \left( {a - {\sum\nolimits_{k
= 1}^{n} {b_{k}}} }  \right)}}}{\prod\limits_{k = 1}^{n} {{\left[
{{\displaystyle\frac{{\Gamma \left( {a - b_{k}}  \right)}}{{\Gamma \left( {a} \right)}}}}
\right]}}}.
\end{equation}

\bigskip
\textbf{Proof}. We carry out the proof by the method mathematical
induction. First, we prove the validity of the equality (\ref{eq888}).

In the case $n=1$ the equality (\ref{eq888}) is obvious.

Let $n = 2$. Since $A(1,2)=A(2,2)=B(1,2)=B(2,2)=m_{2,2}:= i,$
we obtain well-known summation formula (\ref{eq25}):
\begin{equation*}
\label{eq9}
{\sum\limits_{m_{22} = 0}^{\infty}
{{\frac{{\left( {b_{1}}  \right)_{i}}  \left( {b_{1}}  \right)_{i}
}{{\left( {c} \right)_{i}  i!}}}}}:=F\left( {b_{1} ,b_{2} ;a;1} \right)= {\frac{{\Gamma \left( {a -
b_{1} - b_{2}}  \right)\Gamma \left( {a} \right)}}{{\Gamma \left( {a - b_{1}
} \right)\Gamma \left( {a - b_{2}}  \right)}}}.
\end{equation*}

So the formula (\ref{eq888}) works for $n=1$ and $n=2$.

Now we denote the left side of the formula (\ref{eq888}) by

$$T_{n} \left( {a,b_{1} ,...,b_{n}}  \right): =
{\sum\limits_{{\mathop {m_{i,j} = 0}\limits_{(2 \le i \le j \le n)}
}}^{\infty}  {{\displaystyle\frac{{\left( {a} \right)_{A(n,n)}}} {{m_{ij}
!}}}}} {\prod\limits_{k = 1}^{n} {{\displaystyle\frac{{\left( {b_{k}}  \right)_{B(k,n)}
\left( {a - b_{k}}  \right)_{A(k,n) - B(k,n)}}} {{\left( {a}
\right)_{A(k,n)}}} }}} $$

\noindent
and considering fair equality

$$T_{n} \left( {a,b_{1} ,...,b_{n}}  \right) = \Gamma \left( {a -
{\sum\limits_{k = 1}^{n} {b_{k}}} }  \right){\frac{{\Gamma ^{n - 1}\left(
{a} \right)}}{{{\prod\limits_{k = 1}^{n} {\Gamma \left( {a - b_{k}}
\right)}}} }},$$

\noindent
we will prove that
\begin{equation}
\label{eq1001}
T_{n + 1} \left( {a,b_{1} ,...,b_{n + 1}}  \right) = \Gamma \left( {a -
{\sum\limits_{k = 1}^{n + 1} {b_{k}}} }  \right){\frac{{\Gamma ^{n}\left(
{a} \right)}}{{{\prod\limits_{k = 1}^{n + 1} {\Gamma \left( {a - b_{k}}
\right)}}} }}.
\end{equation}

For this aim we will put
\[
T_{n + 1} \left( {a,b_{1} ,...,b_{n + 1}}  \right) =
{\sum\limits_{{\mathop {m_{i,j} = 0}\limits_{(2 \le i \le j \le n+1)}
}}^{\infty}   {{\frac{{\left( {a} \right)_{A(n + 1,n + 1)}
}}{{m_{ij} !}}}}} {\prod\limits_{k = 1}^{n + 1} {{\frac{{\left( {b_{k}}
\right)_{B(k,n + 1)} \left( {a - b_{k}}  \right)_{A(k,n + 1) - B(k,n + 1)}
}}{{\left( {a} \right)_{A(k,n + 1)}}} }}}
\]
and show the validity of the recurrence relation
\begin{equation}
\label{eq11}
T_{n + 1} (a,b_{1} ,...,b_{n + 1} ) = {\prod\limits_{k = 1}^{n} {{\left[
{{\frac{{\Gamma \left( {a} \right)\Gamma \left( {a - b_{k} - b_{n + 1}}
\right)}}{{\Gamma \left( {a - b_{n + 1}}  \right)\Gamma \left( {a - b_{k}}
\right)}}}} \right]}}} \,T_{n} (a - b_{n + 1} ,b_{1} ,...,b_{n} ).
\end{equation}

This process consists of $n$ steps. A detailed look at the first step.

By virtue of the equalities
\[
{\sum\limits_{{\mathop {m_{i,j} = 0}\limits_{(2 \le i \le j \le n+1)}
}}^{\infty}   {}}f(...) = {\sum\limits_{{\mathop {m_{i,j} = 0}\limits_{(2 \le i \le j \le n)}
}}^{\infty}   {{\sum\limits_{{\mathop {m_{i,n+1} = 0}\limits_{(2 \le i \le n+1)}
}}^{\infty}   {}}} }f(...)  = {\sum\limits_{{\mathop {m_{i,j} = 0}\limits_{(2 \le i \le j \le n)}
}}^{\infty}   {{\sum\limits_{{\mathop {m_{i,n+1} = 0}\limits_{(2 \le i \le j \le n)}
}}^{\infty}  }}\sum\limits_{m_{n+1,n+1} = 0}
^{\infty}f(...) }
\]
\noindent
and the properties of functions $A\left( {k,n} \right)$ and $B\left( {k,n}
\right)$ (see formulas (\ref{eq15555}) and (\ref{eq16666})), the right side of equality

$$T_{n + 1} \left( {a,b_{1} ,...,b_{n + 1}}  \right) =
{\sum\limits_{{\mathop {m_{i,j} = 0}\limits_{(2 \le i \le j \le n+1)}
}}^{\infty} {{\frac{{\left( {a} \right)_{A(n + 1,n + 1)}
}}{{m_{ij} !}}}}} {\prod\limits_{k = 1}^{n + 1} {{\frac{{\left( {b_{k}}
\right)_{B(k,n + 1)} \left( {a - b_{k}}  \right)_{A(k,n + 1) - B(k,n + 1)}
}}{{\left( {a} \right)_{A(k,n + 1)}}} }}} $$

\noindent
it is easy to convert to the form
$$
\begin{array}{l}
 T_{n + 1} \left( {a,b_{1} ,...,b_{n + 1}}  \right) =
{\sum\limits_{{\mathop {m_{i,j} = 0}\limits_{(2 \le i \le j \le n)}
}}^{\infty}  {{\displaystyle\frac{{\left( {a - b_{n + 1}}  \right)_{A(n,n)}
\left( {b_{n}}  \right)_{B(n,n)}}} {{m_{ij} !}}}}}  \\
\\
 \cdot{\sum\limits_{{\mathop {m_{i,n+1} = 0}\limits_{(2 \le i  \le n)}
}}^{\infty}  {{\displaystyle\frac{{\left( {b_{n + 1}}  \right)_{m_{2,n + 1} +
... + m_{n,n + 1}}  \left( {a - b_{n}}  \right)_{A(n,n) - B(n,n) + m_{2,n +
1} + ... + m_{n,n + 1}}} } {{m_{i,n + 1} !\left( {a} \right)_{A(n,n) +
m_{2,n + 1} + ... + m_{n,n + 1}}} } }}}  \\
\\
\cdot {\prod\limits_{k = 1}^{n - 1} {{\left[ {{\displaystyle\frac{{\left( {b_{k}}
\right)_{B(k,n) + m_{k + 1,n + 1}}  \left( {a - b_{k}}  \right)_{A(k,n) -
B(k,n) + m_{2,n + 1} + ... + m_{k,n + 1}}} } {{\left( {a} \right)_{A(k,n) +
m_{2,n + 1} + ... + m_{k + 1,n + 1}}} } }S(k,n)} \right]}}} , \\
 \end{array}
$$

\noindent
where
$$
S(k,n) = {\sum\limits_{m_{n + 1,n + 1} = 0}^{\infty}  {{\frac{{\left( {b_{n}
+ B(n,n)} \right)_{m_{n + 1,n + 1}}  \left( {b_{n + 1} + m_{2,n + 1} + ... +
m_{n,n + 1}}  \right)_{m_{n + 1,n + 1}}} } {{m_{n + 1,n + 1} !\left( {a +
A(n,n) + m_{2,n + 1} + ... + m_{n,n + 1}}  \right)_{m_{n + 1,n + 1}}} } }}
}.
$$

It is easy to notice that
$$
S(k,n) = F\left[ b_{n} + B(n,n),b_{n + 1} + m_{2,n + 1} + ... + m_{n,n + 1}
;\right.$$
$$
 \left. a + A(n,n) + m_{2,n + 1} + ... + m_{n,n + 1} ;1 \right].
$$

Applying now the summation formula (\ref{sum}) to the last equality after
elementary transformations we get
$$
 T_{n + 1}^{(1)} \left( {a,b_{1} ,...,b_{n + 1}}  \right) = {\frac{{\Gamma
\left( {a - b_{n} - b_{n + 1}}  \right)\Gamma \left( {a} \right)}}{{\Gamma
\left( {a - b_{n}}  \right)\Gamma \left( {a - b_{n + 1}}
\right)}}}
$$
$$
\cdot{\sum\limits_{{\mathop {m_{i,j} = 0}\limits_{(2 \le i \le j \le n+1)}
}}^{\infty} {{\frac{{\left( {b_{n}}  \right)_{B(n,n)} \left( {a
- b_{n} - b_{n + 1}}  \right)_{A(n,n) - B(n,n)}}} {{m_{ij} !}}}}} {\sum\limits_{{\mathop {m_{i,n+1} = 0}\limits_{(2 \le i \le j \le n)}
}}^{\infty}  {{\frac{{\left( {b_{n + 1}}  \right)_{m_{2,n + 1} +
... + m_{n,n + 1}}} } {{m_{i,n + 1} !}}}}}
$$
$$
\cdot {\prod\limits_{k = 1}^{n - 1}
{{\displaystyle\frac{{\left( {b_{k}}  \right)_{B(k,n) + m_{k + 1,n + 1}}  \left( {a -
b_{k}}  \right)_{A(k,n) - B(k,n) + m_{2,n + 1} + ... + m_{k,n + 1}}
}}{{\left( {a} \right)_{A(k,n) + m_{2,n + 1} + ... + m_{k + 1,n + 1}}} } }}
}.
$$

For definiteness, we denoted the result of the first step of the process
under consideration by $T_{n + 1}^{(1)} \left( {a,b_{1} ,...,b_{n + 1}}
\right)$. We continue the process of proving the recurrence relation (\ref{eq11}).
In each next step, having consistently repeated the reasoning carried out in
the first step, we get
$$
 T_{n + 1}^{(s)} \left( {a,b_{1} ,...,b_{n + 1}}  \right) = {\frac{{\Gamma
^{s}\left( {a} \right)}}{{\Gamma ^{s}\left( {a - b_{n + 1}}
\right)}}}{\prod\limits_{k = n - s + 1}^{n} {{\frac{{\Gamma \left( {a -
b_{k} - b_{n + 1}}  \right)}}{{\Gamma \left( {a - b_{k}}  \right)}}}}}
 $$
 $$
\cdot {\sum\limits_{{\mathop {m_{i,j} = 0}\limits_{(2 \le i \le j \le n)}
}}^{\infty} {{\frac{{1}}{{m_{ij} !}}}{\prod\limits_{k = n - s +
1}^{n} {{\left[ {{\frac{{\left( {b_{k}}  \right)_{B(k,n)} \left( {a - b_{k}
- b_{n + 1}}  \right)_{A(k,n) - B(k,n)}}} {{\left( {a - b_{n + 1}}
\right)_{A(k,n)}}} }} \right]}}}} }
$$
$$
 \cdot {\sum\limits_{{\mathop {m_{i,n+1} = 0}\limits_{(2 \le i  \le n-s+1)}
}}^{\infty} {}} {\frac{{\left( {a - b_{n + 1}}  \right)_{N(n,n)}
\left( {b_{n + 1}}  \right)_{m_{2,n + 1} + ... + m_{n - s + 1,n + 1}}
}}{{m_{ij} !}}}
$$
$$
 \cdot {\prod\limits_{k = 1}^{n - s} {{\left[ {{\frac{{\left( {b_{k}}
\right)_{B(k,n) + m_{k + 1,n + 1}}  \left( {a - b_{k}}  \right)_{A(k,n) -
B(k,n) + m_{2,n + 1} + ... + m_{k,n + 1}}} } {{\left( {a} \right)_{A(k,n) +
m_{2,n + 1} + ... + m_{k + 1,n + 1}}} } }} \right]}}}
$$
\noindent
and in the last step
$$
T_{n + 1}^{(n)} \left( {a,b_{1} ,...,b_{n + 1}}  \right) = {\frac{{\Gamma
^{n}\left( {a} \right)}}{{\Gamma ^{n}\left( {a - b_{n + 1}}
\right)}}}{\prod\limits_{k = 1}^{n} {{\left[ {{\frac{{\Gamma \left( {a -
b_{n + 1} - b_{k}}  \right)}}{{\Gamma \left( {a - b_{k}}  \right)}}}}
\right]}}}
$$
$$
\cdot {\sum\limits_{{\mathop {m_{i,j} = 0}\limits_{(2 \le i \le j \le n)}
}}^{\infty}  {{\frac{{\left( {a - b_{n + 1}}  \right)_{A(n,n)}
}}{{m_{ij} !}}}}} {\prod\limits_{k = 1}^{n} {{\left[ {{\frac{{\left( {b_{k}
} \right)_{B(k,n)} \left( {a - b_{n + 1} - b_{k}}  \right)_{A(k,n) - B(k,n)}
}}{{\left( {a - b_{n + 1}}  \right)_{A(k,n)}}} }} \right]}}} , \\
$$
\noindent
that is
$$
T_{n + 1}^{(n)} \left( {a,b_{1} ,...,b_{n + 1}}  \right) = {\frac{{\Gamma
^{n}\left( {a} \right)}}{{\Gamma ^{n}\left( {a - b_{n + 1}}
\right)}}}{\prod\limits_{k = 1}^{n} {{\left[ {{\frac{{\Gamma \left( {a -
b_{n + 1} - b_{k}}  \right)}}{{\Gamma \left( {a - b_{k}}  \right)}}}}
\right]}}} T_{n} \left( {a - b_{n + 1} ,b_{1} ,...,b_{n}}  \right).
$$

Thus, the validity of the ratio (\ref{eq11}) is established. By the induction
hypothesis, from the (\ref{eq11}) follows the equality
$$
T_{n} \left( {a - b_{n + 1} ,b_{1} ,...,b_{n}}  \right) = \Gamma \left( {a -
b_{n + 1} - {\sum\limits_{k = 1}^{n} {b_{k}}} }  \right){\frac{{\Gamma ^{n -
1}\left( {a - b_{n + 1}}  \right)}}{{{\prod\limits_{k = 1}^{n} {\Gamma
\left( {a - b_{n + 1} - b_{k}}  \right)}}} }}.
$$

Substituting the last expression in (\ref{eq11}) we get the equality (\ref{eq1001}). Therefore, the equality (\ref{eq888}) is true.

The  equality  (\ref{eq8888}) is proved similarly as proof of the equality (\ref{eq888}). Q.E.D.

\bigskip

\textbf{Lemma 3}. The following equalities
$$
{\mathop{\lim} \limits_{\mathop{z_{k} \to 0,}\limits_{k =1,...,n}}}\left\{ z_{1}^{ - b_{1}}  ...z_{n}^{ - b_{n}}  F_{A}^{(n)} \left(
{a,b_{1} ,...,b_{n} ;c_{1} ,...,c_{n} ;1 - {\frac{{1}}{{z_{1}}} },...,1 -
{\frac{{1}}{{z_{n}}} }} \right)\right\}
$$
\begin{equation}
\label{eq12222}
= {\frac{{\Gamma \left( {a - {\sum\nolimits_{k = 1}^{n}
{b_{k}}} }  \right)}}{{\Gamma (a)}}}{\prod\limits_{k = 1}^{n} \left[{{\frac{{\Gamma \left( {c_{k}}
\right)}}{{\Gamma \left( {c_{k} - b_{k}}  \right)}}}}\right]}, a > \sum\limits_{k = 1}^{n}{b_{k}}, b_{k}\neq c_{k}, k=\overline{1,n};
 \end{equation}
\noindent
$$
{\mathop{\lim} \limits_{\mathop{z_{k} \to 0,}\limits_{k =1,...,n}}}\left\{z_{1}^{ - b_{1}}  ...z_{n}^{ - b_{n}}  F_{B}^{(n)} \left( {a_{1}
,...,a_{n} ;b_{1} ,...,b_{n} ;c;1 - {\frac{{1}}{{z_{1}}} },...,1 -
{\frac{{1}}{{z_{n}}} }} \right)\right\}
$$
\begin{equation}
\label{eq12111}
 = {\frac{{\Gamma \left( {c}
\right)}}{{\Gamma \left( {c - {\sum\nolimits_{k = 1}^{n} {b_{k}}} }
\right)}}}{\prod\limits_{k = 1}^{n}\left[ {{\frac{{\Gamma \left( {a_{k} - b_{k}}
\right)}}{{\Gamma \left( {a_{k}}  \right)}}}}\right]}, c > \sum\limits_{k = 1}^{n}{b_{k}}, a_{k} \neq b_{k}, k=\overline{1,n}
\end{equation}
are valid.

\textbf{Proof}. By virtue of the decomposition formula (\ref{eq1222}) we obtain
$$
F_{A}^{(n)} \left( {a,b_{1} ,...,b_{n} ;c_{1} ,...,c_{n} ;1 -
{\frac{{1}}{{z_{1}}} },...,1 - {\frac{{1}}{{z_{n}}} }} \right) =
{\sum\limits_{{\mathop {m_{i,j} = 0}\limits_{(2 \le i \le j \le n)}
}}^{\infty}  {{\frac{{(a)_{A(n,n)}}} {{m_{ij}!} }}}}
$$
\begin{equation}
\label{eq1333}
 \cdot {\prod\limits_{k = 1}^{n}\left[ {{\frac{{(b_{k} )_{B(k,n)}}} {{(c_{k}
)_{B(k,n)}}} }\left( {1 - {\frac{{1}}{{z_{k}}} }} \right)^{B(k,n)}F\left( {a
+ A(k,n),b_{k} + B(k,n);c_{k} + B(k,n);1 - {\frac{{1}}{{z_{k}}} }} \right)}\right]
}.
\end{equation}

Applying now the familiar autotransformation formula (\ref{auto}) to each hypergeometric function included in the sum (\ref{eq1333}), we get
\[
F_{A}^{(n)} \left( {a,b_{1} ,...,b_{n} ;c_{1} ,...,c_{n} ;1 -
{\frac{{1}}{{z_{1}}} },...,1 - {\frac{{1}}{{z_{n}}} }} \right) =
z_{1}^{b_{1}}  ...z_{n}^{b_{n}}  {\sum\limits_{{\mathop {m_{i,j} =
0}\limits_{(2 \le i \le j \le n)}}} ^{\infty}  {{\frac{{(a)_{A(n,n)}
}}{{m_{ij}!} }}}}
\]

\[
 \cdot \prod\limits_{k = 1}^{n} \left[{{\frac{{(b_{k} )_{B(k,n)}}} {{(c_{k}
)_{B(k,n)}}} }\left( {z_{k} - 1} \right)^{B(k,n)}F {\left(
{{\begin{array}{*{20}c}
 {c_{k} - a + B(k,n)
- A(k,n),b_{k} + B(k,n)} ; \hfill \\
 {c_{k} + B(k,n)};  \hfill \\
\end{array}} 1 - z_{k}}  \right)}}\right].
\]

Using the parity property of the sum ${ {B(2,n)+B(3,n)+...+B(n,n)}} $ (see formula (\ref{eq1444})), we calculate the limit
$$
{\mathop{\lim} \limits_{\mathop{z_{k} \to 0,}\limits_{k =1,...,n}}} z_{1}^{ - b_{1}}  ...z_{n}^{ - b_{n}}  F_{A}^{(n)} \left(
{a,b_{1} ,...,b_{n} ;c_{1} ,...,c_{n} ;1 - {\frac{{1}}{{z_{1}}} },...,1 -
{\frac{{1}}{{z_{n}}} }} \right)
$$
$$
= {\sum\limits_{{\mathop {m_{i,j} =
0}\limits_{(2 \le i \le j \le n)}}} ^{\infty}  {{\frac{{(a)_{A(n,n)}
}}{{m_{ij}!} }}}} \prod\limits_{k = 1}^{n} \left[{{\frac{{(b_{k} )_{B(k,n)}}} {{(c_{k}
)_{B(k,n)}}} }F {\left(
{{\begin{array}{*{20}c}
 {c_{k} - a + B(k,n)
- A(k,n),b_{k} + B(k,n)} ; \hfill \\
 {c_{k} + B(k,n)};  \hfill \\
\end{array}} 1}  \right)}}\right]
$$
\noindent
and applying the summation formula (\ref{sum}) to the Gauss hypergeometric functions in the last sum, we obtain the equality (\ref{eq12222}).

The  equality  (\ref{eq12111}) is proved similarly as proof of the equality (\ref{eq12222}). Q.E.D.

\section{Applications of new decomposition formulas to the solution of the boundary value problems}

We consider the equation
\begin{equation}
\label{eq2222}
{\sum\limits_{i = 1}^{m} {u_{x_{i} x_{i}}} }   + {\sum\limits_{k = 1}^{n}
{{\frac{{2\alpha _{k}}} {{x_{k}}} }u_{x_{k}}} }   = 0,
\end{equation}
where $m \ge 2,0 < n \le m;\, \alpha _{k} $ are constants with $\,0 <
2\alpha _{k} < 1 \,$ $\left(k=\overline{1,n}\right)$ in the domain $\Omega$ defined by
$$\Omega\subset {\rm{R}}_m^{n+}:=\{(x_1,...,x_m): x_1>0,...,x_n>0\}.$$  We aim at investigating a Holmgren problem for the equation (\ref{eq2222}).

Let $\Omega \subset {\rm R}_{m}^{n +}  $ be a finite simple-connected domain
bounded by planes $x_{1} = 0,...,x_{n} = 0$ and by the $1/2^n$ part of the $m -
$dimensional sphere $S:$ $x_1^2+...+x_m^2=a^2$. We introduce the notation:
\[
\,\tilde {x}_{k} : = \left( {x_{1} ,...,x_{k - 1} ,x_{k + 1} ,...,x_{n}
,...,x_{m}}  \right) \in S_{k} \subset {\rm R}_{m - 1}^{(n - 1) +}  \subset
{\rm R}^{m - 1} \,\left(k=\overline{1,n}\right).
\]

\textbf{Holmgren problem.} To find a function $u\left( {{x}} \right) \in
C\left( {\bar {\Omega}}  \right) \cap C^{2}\left( {\Omega}  \right)$,
satisfying equation (\ref{eq2222}) in $\Omega $ and conditions
\begin{equation}
\label{eq17}
{\left. {\left( {x_{k}^{2\alpha _{k}}  {\frac{{\partial u}}{{\partial x_{k}
}}}} \right)} \right|}_{x_{k} = 0} = \nu _{k} \left( {\tilde {x}_{k}}
\right),
\,
\,\tilde {x}_{k} \in S_{k} \,\left(k=\overline{1,n}\right),
\end{equation}

\begin{equation}
\label{eq18}
{\left. {u} \right|}_{S} = \varphi \left( {{ x}} \right),
\quad
\,{x} \in \bar {S},
\end{equation}
where $\nu _{k} \left( {\tilde {x}_{k}}  \right)$ and $\varphi \left( {{ x}}
\right)$ are given functions, and, moreover, $\nu _{k} \left( {\tilde
{x}_{k}}  \right)$ can reduce to an infinity of the order less than $1 -2\alpha_1-...-
2\alpha _{n} $ on the boundaries of $S_{k}\,$ $\left(k=\overline{1,n}\right)$.

We find a solution of considered problem using  Green's functions method \cite{A27}.

The Green's function can be represented as
\begin{equation}
\label{eq2144}
G_0\left( {{ x};{ \xi}}  \right) = q_{0} \left( {{ x};{ \xi}}
\right) + q_{0}^{ *}  \left( {{x};{ \xi}}  \right),
\end{equation}
\noindent
where $q_{0} \left( {{x};{\xi}}  \right)$ is the fundamental
solution of equation (\ref{eq2222}), defined by \cite{A288}
\begin{equation*}
q_{0} \left( {{x};{\xi}}  \right) = \gamma _{0} \,r^{ -
2\alpha _{0}} F_{A}^{\left( {n} \right)} \left( {\alpha _{0} ,\alpha _{1}
,...,\alpha _{n} ;2\alpha _{1} ,...,2\alpha _{n} ;\sigma}  \right),
\end{equation*}
where
\[
{x}: = \left( {x_{1} ,...,x_{m}}  \right),
{\xi} : = \left( {\xi _{1} ,...,\xi _{m}}  \right),
\,
{\sigma} : = \left( {\sigma _{1} ,...,\sigma _{n}}  \right);
\]
\begin{equation}
\label{eq16667}
\alpha _{0} = {\frac{{m - 2}}{{2}}} + \alpha _{1} + ... + \alpha _{n} ;\,\,\gamma _{0} = 2^{2\alpha _{0} - m}{\frac{{\Gamma \left( {\alpha _{0}}
\right)}}{{\pi ^{m / 2}}}}{\prod\limits_{k = 1}^{n} {{\frac{{\Gamma \left(
{\alpha _{k}}  \right)}}{{\Gamma \left( {2\alpha _{k}}  \right)}}}}} ,
\end{equation}
\[
r^{2} = {\sum\limits_{i = 1}^{m} {\left( {x_{i} - \xi _{i}}  \right)^{2}}} ,
\,r_{k}^{2} = \left( {x_{k} + \xi _{k}}  \right)^{2} + {\sum\limits_{i = 1,i
\ne k}^{m} {\left( {x_{i} - \xi _{i}}  \right)^{2}}} ,\,
\sigma _{k} = 1 - {\frac{{r_{k}^{2}}} {{r^{2}}}} \, \left(k=\overline{1,n}\right),
\]
a function
\[
q_{0}^{ *}  \left( {{x};{\xi}}  \right) = - \left(
{{\frac{{a}}{{R_{0}}} }} \right)^{2{\alpha}_{0}} q_{0} \left( {{
x};{\bar {\xi}} } \right)
\]
is a regular solution of equation (\ref{eq2222}) in the domain $\Omega $. Here
\[
{ \bar {\xi}} : = \left( {\bar {\xi} _{1} ,...,\bar {\xi} _{m}}  \right),
\bar {\xi} _{i} = {\frac{{a^{2}}}{{R_{0}^{2}}} }\xi _{i} \, \left(i=\overline{1,m}\right);
R_{0}^{2} = \xi _{1}^{2} + ... + \xi _{m}^{2} .
\]

Excise a small ball with its center at ${ \xi} $ and with radius $\rho >
0$ from the domain $\Omega $. Designate the sphere of the excised ball as
$C_{\rho}  $ and by $\Omega _{\rho}  $ denote the remaining part of $\Omega
$.

In deriving an explicit formula for solving the Holmgren problem, the calculation of the following integral plays an important role:
\[
{\int_{C_{\rho}}   {{ x}^{\left( {2\alpha}  \right)}{\left[ {u\left(
{{x}} \right){\frac{{\partial G_0\left( {{x};{\xi}}
\right)}}{{\partial {{\bf n}}}}} - G_0\left( {{ x};{ \xi}}
\right){\frac{{\partial u\left( {{ x}} \right)}}{{\partial {{\bf
n}}}}}} \right]}dC_{\rho}} }
\]
\begin{equation}
\label{eq22}
 = - {\sum\limits_{k = 1}^{n} {{\int_{S_{k}}  {G_0^{\ast}\left(\tilde{x_k}\right)\nu _{k} \left( {\tilde {x}_{k}}
\right)dS_{k}}}  \,}}  + {\int_{S} {{x}^{\left( {2\alpha}
\right)}{\frac{{\partial G_0\left( {{x};{\xi}}
\right)}}{{\partial { {\bf n}}}}}\varphi \left( {\vartheta}
\right)d\vartheta}}
\end{equation}
\noindent
where

\[
{x}^{\left( {2\alpha}  \right)}: = x_{1}^{2\alpha _{1}}
...x_{n}^{2\alpha _{n}}  ,
\,\tilde {x}_{k}^{\left( {2\alpha}  \right)} : = x_{1}^{2\alpha _{1}}  ...x_{k
- 1}^{2\alpha _{k - 1}}  x_{k + 1}^{2\alpha _{k + 1}}  ...x_{n}^{2\alpha
_{n}},
\]
\[
G_0^{\ast}\left(\tilde{x}_k\right):=\tilde {x}_{k}^{\left(
{2\alpha}  \right)} G_0 \left( {x_{1} ,...,x_{k - 1} ,0,x_{k + 1} ,...,x_{m}
;{ \xi}}  \right) \, \left(k=\overline{1,n}\right),
\]
${\bf n}$ is outer normal to $\partial \Omega.$

Since we want to show the application of Lemmas 1-3, therefore, without giving in to details, we discuss only the computation of the following integral
$$
I_{11} = 2\alpha _{0} \gamma _{0} \,\rho ^{ - 2\alpha _{1} - ... - 2\alpha
_{n}} {\int\limits_{0}^{2\pi}  {d\varphi _{m - 1}}} {\int\limits_{0}^{\pi}
{\sin \varphi _{m - 2} d\varphi _{m - 2}}}...
$$
\begin{equation}
\label{l1111} ...{\int\limits_{0}^{\pi}  {
{u\left(\xi_1+\rho\Phi_1,...,\xi_m+\rho\Phi_m\right)}\prod\limits_{i=1}^n\left[\left(\xi_i+\rho\Phi_i\right)^{2\alpha_i}\right]\,F_{A}^{\left( {n} \right)} \left[ {{\sigma}(\rho)}  \right]\sin ^{m
- 2}\varphi _{1} d\varphi _{1}}}  ,
\end{equation}
\noindent
where
\[
\Phi_1=\cos \varphi _{1},\,\Phi_2=\sin\varphi _{1} \cos \varphi _{2},\,\Phi_3=\sin \varphi _{1} \sin \varphi _{2} \cos
\varphi _{3},...,\,\]
\[\Phi_{m-1}=\sin \varphi _{1} \sin \varphi _{2} ...\sin
\varphi _{m - 2} \cos \varphi _{m - 1} ,\,\Phi_m=\sin \varphi _{1} \sin \varphi
_{2} ...\sin \varphi _{m - 2} \sin \varphi _{m - 1};
\]
\[
F_{A}^{\left( {n} \right)} \left( {\sigma _{1\rho}  ,...,\sigma _{n\rho}} \right): =
F_{A}^{\left( {n} \right)} \left( \alpha _{0} + 1,\alpha _{1} ,...,\alpha
_{n} ;2\alpha _{1} ,...,2\alpha _{n} ;{\sigma _{1\rho}  ,...,\sigma _{n\rho}} \right);
\]
\[
r_{k}^{2} = \left( {x_{k} + \xi _{k}}  \right)^{2} + {\sum\limits_{i = 1,i
\ne k}^{m} {\left( {x_{i} - \xi _{i}}  \right)^{2}}} ,
\,
\sigma _{k\rho}  = 1 - {\frac{{r_{k\rho} ^{2}}} {{\rho ^{2}}}} \,\left(k=\overline{1,n}\right).
\]

First we evaluate $F_{A}^{\left( {n} \right)} \left( {\sigma _{1\rho}  ,...,\sigma _{n\rho}}
\right)$. For this aim we use decomposition formula (\ref{eq1222}) and then
formula (\ref{auto}):
\[
F_{A}^{\left( {n} \right)} \left( {\sigma _{1\rho}  ,...,\sigma _{n\rho}} \right) =
{\sum\limits_{{\mathop {m_{i,j} = 0}\limits_{(2 \le i \le j \le n)}
}}^{\infty}  {{\frac{{(\alpha _{0} + 1)_{A(n,n)}}} {{m_{ij} !}}}}
}{\prod\limits_{k = 1}^{n} {{\left[ {{\frac{{(\alpha _{k} )_{B(k,n)}
}}{{(2\alpha _{k} )_{B(k,n)}}} }\left( {1 - {\frac{{r_{k\rho} ^{2}}} {{\rho
^{2}}}}} \right)^{B(k,n)}\left( {{\frac{{r_{k\rho} ^{2}}} {{\rho ^{2}}}}}
\right)^{ - \alpha _{k} - B(k,n)}} \right]}}}
\]
\[
\times {\prod\limits_{k = 1}^{n} {{\left[ {F\left( {2\alpha _{k} - \alpha
_{0} - 1 + B(k,n) - A(k,n),\alpha _{k} + B(k,n);2\alpha _{k} + B(k,n);1 -
{\frac{{r_{k\rho} ^{2}}} {{\rho ^{2}}}}} \right)} \right]}}} ,
\]
\noindent
where $A(k,n)$ and $B(k,n)$ are expressions defined in (\ref{e4111}).

After the elementary evaluations we find
\[
F_{A}^{\left( {n} \right)} \left( {\sigma _{1\rho}  ,...,\sigma _{n\rho}}\right) = \rho
^{2\alpha _{1} + ... + 2\alpha _{n}} {\prod\limits_{k = 1}^{n} {{\left[
{r_{k\rho} ^{ - 2\alpha _{k}}}   \right]}}}  \cdot \aleph ,
\]
\noindent
where
\[
\aleph : = {\sum\limits_{{\mathop {m_{i,j} = 0}\limits_{(2 \le i \le j \le
n)}}} ^{\infty}  {{\frac{{(\alpha _{0} + 1)_{A(n,n)}}} {{m_{ij} !}}}}
}{\prod\limits_{k = 1}^{n} {{\left[ {{\frac{{(\alpha _{k} )_{B(k,n)}
}}{{(2\alpha _{k} )_{B(k,n)}}} }\left( {{\frac{{\rho ^{2}}}{{r_{k\rho} ^{2}
}}} - 1} \right)^{B(k,n)}} \right]}}}
\]
\[
\times {\prod\limits_{k = 1}^{n} {{\left[ {F\left( {2\alpha _{k} - \alpha
_{0} - 1 + B(k,n) - A(k,n),\alpha _{k} + B(k,n);2\alpha _{k} + B(k,n);1 -
{\frac{{\rho ^{2}}}{{r_{k\rho} ^{2}}} }} \right)} \right]}}} .
\]

It is easy to see that when $\rho \to 0$ the function $\aleph\, $ becomes an
expression that does not depend on ${\rm x}$ and ${\rm \xi} $. Indeed,
taking into account
the parity property of the sum ${ {B(2,n)+B(3,n)+...+B(n,n)}} $ (see formula (\ref{eq1444})), we have
\[
{\mathop {\lim} \limits_{\rho \to 0}} \aleph : = {\sum\limits_{{\mathop
{m_{i,j} = 0}\limits_{(2 \le i \le j \le n)}}} ^{\infty}  {{\frac{{(\alpha
_{0} + 1)_{A(n,n)}}} {{m_{ij} !}}}}} {\prod\limits_{k = 1}^{n} {{\left[
{{\frac{{(\alpha _{k} )_{B(k,n)}}} {{(2\alpha _{k} )_{B(k,n)}}} }} \right]}}
}
\]
\begin{equation}
\label{sum1}
\times {\prod\limits_{k = 1}^{n} {{\left[ {F\left( {2\alpha _{k} - \alpha
_{0} - 1 + B(k,n) - A(k,n),\alpha _{k} + B(k,n);2\alpha _{k} + B(k,n);1}
\right)} \right]}}} .
\end{equation}

Applying now the summation formula (\ref{sum}) to each hypergeometric function
$F\left( {a,b;c;1} \right)$ in the sum (\ref{sum1}), we get
\[
{\mathop {\lim} \limits_{\rho \to 0}} \aleph : = {\frac{{1}}{{\Gamma \left(
{\alpha _{0} + 1} \right)}}}{\sum\limits_{{\mathop {m_{i,j} = 0}\limits_{(2
\le i \le j \le n)}}} ^{\infty}  {{\frac{{\Gamma (\alpha _{0} + 1 +
N(n,n))}}{{m_{ij} !}}}}}
\]
\[
\cdot {\prod\limits_{k = 1}^{n} {{\left[ {{\frac{{\Gamma \left( {2\alpha
_{k}}  \right)\Gamma (\alpha _{k} + M(k,n))\Gamma \left( {\alpha _{0} + 1 -
\alpha _{k} + N(k,n) - M(k,n)} \right)}}{{\Gamma ^{2}\left( {\alpha _{k}}
\right)\Gamma \left( {\alpha _{0} + 1 + N(k,n)} \right)}}}} \right]}}} .
\]

Taking into account the identity (\ref{eq888}) we obtain
\begin{equation}
\label{eq3000}
{\mathop {\lim} \limits_{\rho \to 0}} \aleph = {\frac{{\Gamma \left( {m / 2}
\right)}}{{\Gamma \left( {\alpha _{0} + 1} \right)}}}{\prod\limits_{i =
1}^{n} {{\frac{{\Gamma \left( {2\alpha _{k}}  \right)}}{{\Gamma \left(
{\alpha _{k}}  \right)}}}}} .
\end{equation}

Now we consider an integral
\[
L_{m} = {\int\limits_{0}^{2\pi}  {d\varphi _{m - 1}}}  {\int\limits_{0}^{\pi}
{\sin \varphi _{m - 2} d\varphi _{m - 2}}}  {\int\limits_{0}^{\pi}  {\sin ^{2}\varphi
_{m - 3} d\varphi _{m - 3}}}  ...{\int\limits_{0}^{\pi}  {\sin ^{m - 2}\varphi
_{1} d\varphi _{1}}},
\]
with elementary transformations it is not difficult to establish that
\begin{equation}
\label{eq3111}
L_{2m} = {\frac{{2\,\pi ^{m}}}{{(m - 1)!}}},\,\,L_{2m + 1} = {\frac{{2^{m +
1}\,\pi ^{m}}}{{(2m - 1)!!}}},\,\,\,m = 1,2,3,...
\end{equation}

If we take into account (\ref{eq16667}), (\ref{l1111}), (\ref{eq3000}) and (\ref{eq3111}) , then from (\ref{eq22}) we will have
\[
{\mathop {\lim} \limits_{\rho \to 0}} I_{11} = u\left( {\xi}  \right).
\]

So  we can write the solution of the Holmgren problem as follows:
\begin{equation}
\label{9999}
u\left( {\xi}  \right) = - {\sum\limits_{k = 1}^{n} {{\int_{S_{k}}  { G_0^{\ast}  \left( {\tilde
{x}_{k};{\xi}}  \right)\nu _{k} \left( {\tilde
{x}_{k}}  \right)dS_{k}}}  \,}}
 + {\int_{S} {x^{\left( {2\alpha}
\right)}{\frac{{\partial G_0\left( {{ x};{\xi}}
\right)}}{{\partial {\rm {\bf n}}}}}\varphi \left( {{x}} \right)dS}},
\end{equation}
where
\begin{equation*}
G_0^{\ast}  \left( {\tilde
{x}_{k};{\xi}}  \right)=\gamma _{0}\tilde {x}_{k}^{\left( {2\alpha}  \right)}{\left\{ {{\frac{{F_{A}^{\left( {n - 1} \right)} {\left[
{{\begin{array}{*{20}c}
 {{\alpha} _{0} ,\alpha _{1} ,...,\alpha _{k - 1} ,\alpha
_{k + 1} ,...,\alpha _{n} ;} \hfill \\
 {2\alpha _{1} ,...,2\alpha _{k - 1} ,2\alpha _{k + 1} ,...,
2\alpha _{n} ;} \hfill \\
\end{array}} \sigma _{0}}  \right]}}}{{{\left[ {\xi _{k}^{2} +
{\sum\limits_{i = 1,i \ne k}^{m} {\left( {\xi _{i} - x_{i}}  \right)^{2}}}}
\right]}^{\,{\alpha} _{0}}} }}} \right.}
\end{equation*}
\begin{equation*}
\left.{ - {\frac{{F_{A}^{\left( {n - 1} \right)} {\left[
{{\begin{array}{*{20}c}
 {{\alpha}_{0} ,\alpha _{1} ,...,\alpha _{k - 1} ,\alpha
_{k + 1} ,...,\alpha _{n} ;} \hfill \\
 {2\alpha _{1} ,...,2\alpha _{k - 1} ,2\alpha _{k + 1} ,...,
2\alpha _{n} ;} \hfill \\
\end{array}} \bar {\sigma} _{0}}  \right]}}}{{{\left[\sum\limits_{i = 1,i
\ne k}^{m} \left( a - \displaystyle\frac{x_i\xi_i}{a} \right)^{2}  +
\displaystyle\frac{1}{a^2}{\sum\limits_{i = 1,i \ne k}^{m} {{\sum\limits_{j=1,j
\ne i}^{m} {x_{i}^{2} \xi _{j}^{2}}}}}-(m-2)a^{2}
\right]}^{\, {\alpha} _{0}}}}}}\right\},
\end{equation*}
$G_0\left( {{ x};{\xi}}
\right)$  is the Green's function, defined by (\ref{eq2144}).

In conclusion, we note precisely because of the decomposition formula (\ref{eq1222}), the summation formula (\ref{eq888}) and the limit value (\ref{eq12222}) that we managed to write out the solution of the Holmgren problem with conditions (\ref{eq17}) and (\ref{eq18}) for the equation  (\ref{eq2222}) in an explicit form (\ref{9999}).

 \end{document}